\documentclass{article}
\usepackage{amssymb}
\usepackage{amsmath}
\usepackage{amsfonts}       
\usepackage{color} 
\newcommand{\email}[1]{{\small E-mail: {\textsf {#1}}}}
\newtheorem{theo}{Theorem}[section]
\newtheorem{prop}[theo]{Proposition}
\newtheorem{lem}[theo]{Lemma}
\newtheorem{cor}[theo]{Corollary}

\newtheorem{defi}{Definition}

\newcommand{\mysection}[1]{\section{#1} \setcounter{equation}{0}}

\newcommand{\proof}{{\sc Proof.} \quad}
\newcommand{\proofc}[1]{{\sc Proof {#1}.}}        

\newcommand{\R}{\mathbb{R}}

\newcommand{\be}{\begin{equation} \label}
\newcommand{\ee}{\end{equation}}
\newcommand{\bes}{\begin{equation} \begin{array}{c} \label}
\newcommand{\ees}{\end{array} \end{equation}}
\newcommand{\bea}{\begin{eqnarray}\label}
\newcommand{\eea}{\end{eqnarray}}
\newcommand{\beas}{\begin{eqnarray} \begin{array}{rcl} \label}
\newcommand{\eeas}{\end{array} \end{eqnarray}}
\newcommand{\bas}{\begin{eqnarray*}}\newcommand{\eas}{\end{eqnarray*}}
\newcommand{\bass}{\begin{eqnarray*} \begin{array}{rcl}}
\newcommand{\eass}{\end{array} \end{eqnarray*}}
\newcommand{\basss}{\begin{eqnarray*} \begin{array}{c}}
\newcommand{\easss}{\end{array} \end{eqnarray*}}
\newcommand{\qed}{{}\hfill //// \\[15pt]}
\newcommand{\bit}{\begin{itemize}}
\newcommand{\eit}{\end{itemize}}

\newcommand{\eps}{\varepsilon}
\newcommand{\abs}{\\[3mm]}

\newcommand{\pO}{\partial\Omega}

\begin{document}
\title{Refined asymptotics for the infinite heat equation with homogeneous Dirichlet boundary conditions 
}
\author{Philippe Lauren\c{c}ot\footnote{Institut de Math\'ematiques de Toulouse, CNRS UMR~5219, Universit\'e de Toulouse, F--31062 Toulouse cedex 9, France. \email{laurenco@math.univ-toulouse.fr}} \, \& 
Christian Stinner\footnote{Fakult\"{a}t f\"{u}r Mathematik, Universit\"at Duisburg-Essen, D--45117 Essen, Germany. \email{christian.stinner@uni-due.de}}}
\date{\today}
\maketitle
\begin{abstract}
\noindent  The nonnegative viscosity solutions to the infinite heat equation with homogeneous Dirichlet boundary conditions
  are shown to converge as $t \to \infty$ to a uniquely determined limit after a suitable time rescaling. The proof relies on
  the half-relaxed limits technique as well as interior positivity estimates and boundary estimates. The expansion of the
  support is also studied. 
\end{abstract}
%
%
%
%
%
%
%%%%%%%%%%%%%%%%%%%%%%%%%%%%%%%%%%%%%%%%%%%%%%%%%%%%%%%%%%%%%%%%%%%%%%%%%%%%%%%%%%
%%%%%%%%%%%%%%%%%%%%%%%%%%%%%%%%%%%%%%%%%%%%%%%%%%%%%%%%%%%%%%%%%%%%%%%%%%%%%%%%%%
\mysection{Introduction}
%%%%%%%%%%%%%%%%%%%%%%%%%%%%%%%%%%%%%%%%%%%%%%%%%%%%%%%%%%%%%%%%%%%%%%%%%%%%%%%%%%
%%%%%%%%%%%%%%%%%%%%%%%%%%%%%%%%%%%%%%%%%%%%%%%%%%%%%%%%%%%%%%%%%%%%%%%%%%%%%%%%%%
Since the pioneering work by Aronsson \cite{A67}, the infinity-Laplacian $\Delta_\infty$ defined by
$$\Delta_\infty u := \langle D^2 u \nabla u, \nabla u \rangle = \sum\limits_{i,j = 1}^N \frac{\partial^2 u}{\partial x_i 
  \partial x_j}\frac{\partial u}{\partial x_i} \frac{\partial u}{\partial x_j}$$
has been the subject of several studies, in particular due to its relationship to the theory of absolutely minimizing
Lipschitz extensions \cite{A67, ACJ04, Cr08}. More recently, a parabolic equation involving the infinity-Laplacian (the
infinite heat equation)
\begin{equation}\label{P}
  \partial_t u = \Delta_\infty u,        \qquad (t,x ) \in (0,\infty) \times \Omega,
\end{equation}
has been considered in \cite{AJK09, AS08, CW03}. When $\Omega \subset \mathbb{R}^N$ is a bounded domain and
\eqref{P} is supplemented with nonhomogeneous Dirichlet
boundary conditions, the large time behaviour of solutions to \eqref{P} is investigated in \cite{AJK09} and convergence as $t \to \infty$ to
the unique steady state is shown. Furthermore, for homogeneous Dirichlet boundary conditions
\begin{equation}\label{PBC}
  u = 0,   \qquad (t,x) \in (0,\infty) \times \partial \Omega,
\end{equation}
and nonnegative initial condition
\begin{equation}\label{PIC}
        u(0,x) = u_0(x),         \qquad x\in \bar{\Omega}, 
\end{equation}
satisfying
\begin{equation}\label{ic}
  u_0 \in C_0 (\bar{\Omega}) := \{ f \in C(\bar{\Omega}) \, : \; f=0 \mbox{ on } \partial \Omega \}, \; u_0 \ge 0, \; u_0 \not\equiv 0,
\end{equation}
a precise temporal decay rate is given for the $L^\infty$-norm of $u$, namely
\begin{equation}\label{1.1.0}
  C_1^{-1}(t+1)^{-1/2} \le \| u(t, \cdot) \|_{L^\infty (\Omega)} \le C_1(t+1)^{-1/2} \quad\mbox{for all } t>0
\end{equation}
with some $C_1 \ge 1$ depending on $u_0$ and $\Omega$, the unique steady state of \eqref{P}-\eqref{PBC} being zero in that
case.

The purpose of this note is to improve \eqref{1.1.0} by identifying the limit of $t^{1/2} u(t, \cdot)$ as $t \to \infty$
(see Theorem~\ref{theo1.2} below). We also provide additional information on the propagation of the positivity set of $u$
as time goes by. 

Before stating our main result we first recall that the infinity-Laplacian is a quasilinear and degenerate elliptic
operator which is not in divergence form and a suitable framework to study the well-posedness of the infinite heat equation
is the theory of viscosity solutions (see e.g. \cite{CIL92}). Within this framework the well-posedness of \eqref{P}-\eqref{PIC}
has been established in \cite{AS08} when $\Omega$ fulfills the \textit{uniform exterior sphere condition}: 
\begin{equation}\label{esc}
  \begin{array}{l}
  \mbox{For all } x_0 \in \partial \Omega \mbox{ there exists } y_0 \in \mathbb{R}^N \mbox{ such that } |x_0 - y_0| = R 
  \mbox{ and } \\ \{x \in \mathbb{R}^N \, : \, |x-y_0| <R \} \cap \Omega = \emptyset \mbox{ for some positive constant } R
  \mbox{ independent of } x_0. \end{array}
\end{equation}

Introducing
\begin{equation}\label{F}
  F(s,p,X) := s - \langle X p,p \rangle \quad\mbox{for } s \in \mathbb{R}, p \in \mathbb{R}^N, X \in \mathcal{S} (N),
\end{equation} 
where $\mathcal{S}(N)$ denotes the set of all symmetric $N \times N$ matrices, the definition of viscosity solutions to
\eqref{P}-\eqref{PIC} reads \cite{AJK09,AS08}:
\begin{defi}\label{defi1.1}
  Let $Q := (0,\infty) \times \Omega \subset \mathbb{R}^{N+1}$ and let $USC(\bar{Q})$ and $LSC(\bar{Q})$ denote the set of
  upper semicontinuous and lower semicontinuous functions from $\bar{Q}$ into $\mathbb{R}$, respectively. A function
  $u \in USC(\bar{Q})$ is a viscosity subsolution to \eqref{P}-\eqref{PIC} in $Q$ if 
  \begin{itemize}
    \item[(a)] $F(s,p,X) \le 0$ is satisfied for all $(s,p,X) \in \mathcal{P}^{2,+} u(t_0,x_0)$ and all $(t_0,x_0) \in Q$,
      where
      \begin{eqnarray*}
       \mathcal{P}^{2,+} u(t_0,x_0) &:=& \Big\{ (s,p,X) \in \mathbb{R} \times \mathbb{R}^N \times \mathcal{S}(N) \; : \; 
       u(t,x) \le u(t_0,x_0) + s(t-t_0) \\
      & & + \langle p, x-x_0 \rangle + \frac{1}{2} \langle X (x-x_0), x-x_0 \rangle + o(|t-t_0| + |x-x_0|^2) \\
      & & \mbox{ as } (t,x) \to (t_0,x_0) \Big\},
      \end{eqnarray*} 
    \item[(b)] $u \le 0$ on $(0,\infty) \times \partial \Omega$, 
    \item[(c)] $u (0,x) \le u_0(x)$ for $x \in \bar{\Omega}$. 
  \end{itemize}
  Similarly, $u \in LSC(\bar{Q})$ is a viscosity supersolution to \eqref{P}-\eqref{PIC} in $Q$ if 
  $F(s,p,X) \ge 0$ for all $(s,p,X) \in \mathcal{P}^{2,-} u(t_0,x_0) := - \mathcal{P}^{2,+} (-u)(t_0,x_0)$ and $(t_0,x_0) \in 
  Q$, $u \ge 0$ on $(0,\infty) \times \partial \Omega$ and $u (0,x) \ge u_0(x)$ for $x \in \bar{\Omega}$.
  
  Finally, $u \in C(\bar{Q})$ is a viscosity solution to \eqref{P}-\eqref{PIC} if it is a viscosity subsolution and a 
  viscosity supersolution to \eqref{P}-\eqref{PIC}.
\end{defi}
\bigskip
With this definition, the well-posedness of \eqref{P}-\eqref{PIC} is shown in \cite[Theorems~2.3 and 2.5]{AS08} and the
asymptotic behaviour of nonnegative solutions is obtained in \cite[Theorem~5]{AJK09}. We gather these results in the next
theorem.
\begin{theo}\label{theo1.1} (\cite{AJK09,AS08})
  Let $\Omega \subset \mathbb{R}^N$ be a bounded domain such that \eqref{esc} is satisfied and assume \eqref{ic}. Then there
  is a unique nonnegative viscosity solution $u$ to \eqref{P}-\eqref{PIC}. Moreover, $u (t, \cdot)$ converges to zero as $t \to \infty$
  in the sense that there exists a constant $C_1 \ge 1$ such that
  \begin{equation}\label{1.1.1}
    C_1^{-1}(t+1)^{-1/2} \le \| u(t, \cdot) \|_{L^\infty (\Omega)} \le C_1(t+1)^{-1/2} \quad\mbox{for all } t>0.
  \end{equation}
\end{theo}
\bigskip
Our improvement of \eqref{1.1.1} then reads:
\begin{theo}\label{theo1.2}
  Suppose $\Omega \subset \mathbb{R}^N$ is a bounded domain fulfilling \eqref{esc} and assume \eqref{ic}. If
  $u$ denotes the viscosity solution to \eqref{P}-\eqref{PIC}, then 
  \begin{equation}\label{1.2.1}
    \lim\limits_{t \to \infty} \| t^{1/2} u(t, \cdot) - f_\infty \|_{L^\infty (\Omega)} = 0,
  \end{equation}
  where $f_\infty$ is the unique positive viscosity solution to
  \begin{equation}\label{1.2.2}
    - \Delta_\infty f_\infty - \frac{f_\infty}{2} = 0 \mbox{ in } \Omega, \quad f_\infty>0 \mbox{ in } \Omega, \quad 
    f_\infty = 0 \mbox{ on } \partial \Omega.
  \end{equation}
\end{theo}
\bigskip
Theorem~\ref{theo1.2} not only gives the convergence of $t^{1/2} u(t, \cdot)$ as $t \to \infty$, but also provides the
existence and uniqueness of the positive solution $f_\infty$ to \eqref{1.2.2} in $C_0 (\bar{\Omega})$. An interesting
consequence of
\eqref{1.2.2} is that the function $(t,x) \mapsto t^{-1/2} f_\infty (x)$ is a separate variables solution to 
\eqref{P}-\eqref{PBC} with an initial data being identically infinite in $\Omega$. Similar solutions are already known to
exist for other parabolic equations such as the porous medium equation $\partial_t u = \Delta u^m$, $m>1$, or the 
$p$-Laplacian equation $\partial_t u = {\rm div} (|\nabla u|^{p-2} \nabla u)$, $p > 2$, (see 
\cite{AP81, CLS89, DK88, MV94, Va07} for instance). They play an important role in the description of the large time dynamics 
\cite{AP81, Ga04, Va07} and also provide
universal bounds (and are thus called \textit{friendly giants}). The function $(t,x) \mapsto t^{-1/2} f_\infty (x)$ is a
friendly giant for the infinite heat equation \eqref{P}-\eqref{PIC} and we have the following universal bound.
\begin{cor}\label{cor1.3}
  Suppose $\Omega \subset \mathbb{R}^N$ is a bounded domain fulfilling \eqref{esc} and assume \eqref{ic}. If
  $u$ denotes the viscosity solution to \eqref{P}-\eqref{PIC}, then
  \begin{equation}\label{1.3.1}
    u(t,x) \le t^{-1/2} f_\infty (x) \qquad\mbox{for } (t,x) \in (0,\infty) \times \bar{\Omega},
  \end{equation}
  the function $f_\infty$ being defined in Theorem~\ref{theo1.2}.
\end{cor}
\bigskip
The proof of Theorem~\ref{theo1.2} and Corollary~\ref{cor1.3} involves several steps: According to \eqref{1.1.1} the 
evolution of $u(t, \cdot)$ takes place on a time scale of order $t^{-1/2}$ and we first introduce a rescaled version $v$ of
$u$ defined by $u(t,x) = t^{-1/2} v ( \ln t,x)$. The outcome of Theorem~\ref{theo1.2} is then the convergence of $v(s,\cdot)$ 
to the time-independent function $f_\infty$ as $s \to \infty$. To establish such a convergence, we use the half-relaxed 
limits technique introduced in \cite{BP88} which is well-suited here as we have rather scarce information on $v (s,\cdot)$ as
$s \to \infty$. This requires however a strong comparison principle for the limit problem \eqref{1.2.2} which will be
established in Section~\ref{ufg}, under an additional positivity assumption, and furthermore implies the uniqueness of 
$f_\infty$. That the half-relaxed limits indeed enjoy this positivity property has to be proved as a preliminary step
and follows from the observation that $v(s, \cdot)$ is non-decreasing with time and eventually becomes positive in $\Omega$ 
(see Section~\ref{ptm}). At this point, boundary estimates are also needed to ensure that the half-relaxed limits vanish
on $\partial \Omega$ and are shown by constructing suitable barrier functions. Thanks to these results, we deduce that 
the half-relaxed limits coincide, which implies that $v(s, \cdot)$ converges as $s \to \infty$ and the existence of a
positive solution $f_\infty$ to \eqref{1.2.2} as well (see Section~\ref{cv}). We emphasize here that the existence of a positive 
solution to \eqref{1.2.2} is a consequence of the dynamical properties of $v$ and was seemingly not known previously. Finally, 
Corollary~\ref{cor1.3} is a consequence of Theorem~\ref{theo1.2} and the time monotonicity of $v$ (see Section~\ref{cv}).

Additionally, in Section~\ref{app} we investigate further positivity properties of the solution $u$ to \eqref{P}-\eqref{PIC}.
We show that $u(t,\cdot)$ becomes positive in $\Omega$ after a finite time if $\Omega$ satisfies an additional uniform interior
sphere condition. Aside from this, $u$ may have a positive waiting time if the initial data are flat on the boundary of their
support, namely the support of $u(t, \cdot)$ will be equal to that of $u_0$ for small times.    

\medskip

For further use, we introduce the following notation:
Given $x \in \bar{\Omega}$, let $d(x, \partial \Omega) := dist(x, \partial \Omega)$ denote the distance to the boundary. 
Moreover, for $x \in \mathbb{R}^N$ and $r>0$ we define $B(x,r) := \{ y \in \mathbb{R}^N \, : \, |y-x| < r\}$ to be the ball of 
radius $r$ centered at $x$.
%
%
%
%
%
%
%
%
%%%%%%%%%%%%%%%%%%%%%%%%%%%%%%%%%%%%%%%%%%%%%%%%%%%%%%%%%%%%%%%%%%%%%%%%%%%%%%%%%%
%%%%%%%%%%%%%%%%%%%%%%%%%%%%%%%%%%%%%%%%%%%%%%%%%%%%%%%%%%%%%%%%%%%%%%%%%%%%%%%%%%
\mysection{Uniqueness of the friendly giant}\label{ufg}
%%%%%%%%%%%%%%%%%%%%%%%%%%%%%%%%%%%%%%%%%%%%%%%%%%%%%%%%%%%%%%%%%%%%%%%%%%%%%%%%%%
%%%%%%%%%%%%%%%%%%%%%%%%%%%%%%%%%%%%%%%%%%%%%%%%%%%%%%%%%%%%%%%%%%%%%%%%%%%%%%%%%%
In this section we show that the friendly giant is unique. This will be a consequence of the following more general comparison
lemma.
\begin{lem}\label{lem2.1}
  Let $w \in USC( \bar{\Omega})$ and $W \in LSC(\bar{\Omega})$ be respectively a bounded viscosity subsolution and a bounded
  viscosity supersolution to
  \begin{equation}\label{2.1.1}
    - \Delta_\infty \zeta - \frac{\zeta}{2} = 0 \quad\mbox{in } \Omega
  \end{equation}
  such that
  \begin{eqnarray}
    & & w(x) = W(x) = 0 \quad\mbox{for } x \in \partial \Omega, \label{2.1.2} \\
    & & W(x) > 0 \hspace*{+12mm} \quad\mbox{for } x \in \Omega. \label{2.1.3}
  \end{eqnarray}
  Then 
  \begin{equation}\label{2.1.4}
    w \le W \quad\mbox{in } \Omega.
  \end{equation}
\end{lem}
\proof We fix $N_0 \in \mathbb{N}$ large enough such that $\Omega_n := \{ x \in \Omega \, : \, d(x, \pO) > 1/n \}$ is a
  nonempty open subset of $\Omega$ for all integer $n \ge N_0$. Let $n \ge N_0$. Since $\bar{\Omega}_n$ is
  compact and $W \in LSC(\bar{\Omega})$, $W$ has a minimum in $\bar{\Omega}_n$ and the positivity of $W$
  in $\bar{\Omega}_n$ implies that
  \begin{equation}\label{2.1.5}
    \mu_n := \min\limits_{\bar{\Omega}_n} \, W >0.
  \end{equation}
  Similarly, the compactness of $\bar{\Omega} \setminus \Omega_n$ and the upper semicontinuity and boundedness of $w$ 
  ensure that $w$ has a point of maximum $x_n$ in $\bar{\Omega} \setminus \Omega_n$ and we set
  \begin{equation}\label{2.1.6}
    \eta_n := \max\limits_{\bar{\Omega} \setminus \Omega_n} \, w = w(x_n) \ge 0,
  \end{equation}
  the nonnegativity of $\eta_n$ being a consequence of the fact that $w$ vanishes of $w$ on $\partial \Omega$. We next claim
  that
  \begin{equation}\label{2.1.7}
    \lim\limits_{n \to \infty} \eta_n = 0.
  \end{equation}
  Indeed, owing to the compactness of $\bar{\Omega}$ and the definition of $\Omega_n$ there are $y \in \partial \Omega$
  and a subsequence of $(x_n)_{n \in \mathbb{N}}$ (not relabeled) such that $x_n \to y$ as $n \to \infty$. Since $w(y) = 0$,
  we deduce from the upper semicontinuity of $w$ that
  $$\limsup_{x \to y} w(x) = \lim\limits_{\eps \searrow 0} \; \sup \{ w(x) \, : \, x \in B(y, \eps) \cap \bar{\Omega} \} \le 
  0.$$
  Given $\eps >0$, there is $n_\eps$ such that $x_n \in B(y, \eps) \cap \bar{\Omega}$ for all $n \ge n_\eps$. Hence,
  $$\limsup\limits_{n \to \infty} \eta_n \le \sup \{ w(x) \, : \, x \in B(y, \eps) \cap \bar{\Omega} \}$$
  and letting $\eps \searrow 0$ and using \eqref{2.1.6} allow us to conclude that
  $$0 \le \limsup\limits_{n \to \infty} \eta_n \le 0.$$
  This shows that a subsequence of $(\eta_n)_{n \ge N_0}$ converges to zero and the claim \eqref{2.1.7} follows by
  noticing that $(\eta_n)_{n \ge N_0}$ is a nonincreasing sequence.  
  
  Next, fix $s \in (0,\infty)$. For $\delta >0$ and $n \ge N_0$, we define
  \begin{eqnarray*}
    && \begin{array}{ll}
    z_n (t,x) := (t+s)^{-1/2} w(x) - s^{-1/2} \eta_n, & \quad (t,x) \in [0,\infty) \times \bar{\Omega}, \\
    Z_\delta (t,x) := (t+ \delta )^{-1/2} W(x), & \quad (t,x) \in [0,\infty) \times \bar{\Omega}. \end{array}
  \end{eqnarray*}
  Then $z_n$ and $Z_\delta$ are respectively a bounded usc viscosity subsolution and a bounded lsc viscosity supersolution
  to \eqref{P} with
  $$Z_\delta (t,x) = 0 \ge - s^{-1/2} \eta_n = z_n(t,x), \quad (t,x) \in (0,\infty) \times \partial \Omega.$$
  In addition, if
  \begin{equation}\label{2.1.8}
    0 < \delta < \left( \frac{\mu_n}{1 + \| w \|_{L^\infty (\Omega)}} \right)^2 s
  \end{equation}  
  we have 
  $$Z_\delta (0,x) = \delta^{-1/2} W(x) \ge \delta^{-1/2} \mu_n \ge s^{-1/2} \| w \|_{L^\infty (\Omega)} \ge z_n (0,x)
  \quad\mbox{for } x \in \Omega_n$$
  and
  $$Z_\delta (0,x) \ge 0 \ge s^{-1/2} (w(x) - \eta_n) = z_n (0,x) \quad\mbox{for } x \in \bar{\Omega} \setminus \Omega_n.$$
  We are then in a position to apply the comparison principle \cite[Theorem~8.2]{CIL92} to deduce that
  \begin{equation}\label{2.1.9}
    z_n(t,x) \le Z_\delta (t,x), \quad (t,x) \in [0,\infty) \times \bar{\Omega}, 
  \end{equation}
  for any $\delta >0$ and $n \ge N_0$ satisfying \eqref{2.1.8}. According to \eqref{2.1.8}, the parameter $\delta$ can be 
  taken arbitrarily small and we deduce from \eqref{2.1.9} that 
  $$(t+s)^{-1/2} w(x) - s^{-1/2} \eta_n \le t^{-1/2} W(x), \quad (t,x) \in (0,\infty) \times \bar{\Omega},$$
  for $n \ge N_0$. We next pass to the limit as $n \to \infty$ with the help of \eqref{2.1.7} to conclude that
  $$(t+s)^{-1/2} w(x) \le t^{-1/2} W(x), \quad (t,x) \in (0,\infty) \times \bar{\Omega}.$$  
  Finally, as $s>0$ is arbitrary, we may let $s \searrow 0$ and take $t=1$ in the above inequality to complete the proof.
\qed
Now the uniqueness of the friendly giant is a straightforward consequence of Lemma~\ref{lem2.1}.
\begin{cor}\label{cor2.2}
  There is at most one positive viscosity solution to \eqref{1.2.2} in $C_0 (\bar{\Omega})$.
\end{cor}
\bigskip
%
%
%
%
%%%%%%%%%%%%%%%%%%%%%%%%%%%%%%%%%%%%%%%%%%%%%%%%%%%%%%%%%%%%%%%%%%%%%%%%%%%%%%%%%%
%%%%%%%%%%%%%%%%%%%%%%%%%%%%%%%%%%%%%%%%%%%%%%%%%%%%%%%%%%%%%%%%%%%%%%%%%%%%%%%%%%
\mysection{Large time behaviour}\label{ltb}
%%%%%%%%%%%%%%%%%%%%%%%%%%%%%%%%%%%%%%%%%%%%%%%%%%%%%%%%%%%%%%%%%%%%%%%%%%%%%%%%%%
%%%%%%%%%%%%%%%%%%%%%%%%%%%%%%%%%%%%%%%%%%%%%%%%%%%%%%%%%%%%%%%%%%%%%%%%%%%%%%%%%%
%
In this section, we assume that $\Omega$ is a bounded domain fulfilling \eqref{esc} and that $u_0$ satisfies \eqref{ic}.
Let $u$ be the corresponding viscosity solution to \eqref{P}-\eqref{PIC}.
In order to investigate the asymptotic behaviour of $u$ as stated in Theorem~\ref{theo1.2} we
introduce the scaling variable $s= \ln t$, $t > 0$, and the rescaled unknown function $v$ defined by
\begin{equation}\label{3.1}
  u(t,x) = t^{-1/2} v ( \ln t, x ), \quad (t,x) \in (0,\infty) \times \bar{\Omega}.
\end{equation}
It is easy to check that $v$ is the viscosity solution to
\begin{eqnarray}
        \partial_s v & = & \Delta_\infty v - \frac{v}{2},   \hspace*{+8mm}\qquad (s,x ) \in (0,\infty) \times \Omega, 
        \label{v} \\
        v & = & 0,         \hspace*{+21mm} \qquad (s,x) \in (0,\infty) \times \partial \Omega, \label{vbc}\\
        v(0,x) & = & v_0(x) := u(1,x),       \qquad x\in \bar{\Omega}, \label{vic} 
\end{eqnarray}
while it readily follows from \eqref{1.1.1} and \eqref{3.1} that
\begin{equation}\label{3.1a}
  0 \le v(s,x) \le C_1, \qquad (s,x) \in [0,\infty) \times \bar{\Omega}.
\end{equation}
%
%%%%%%%%%%%%%%%%%%%%%%%%%%%%%%%%%%%%%%%%%%%%%%%%%%%%%%%%%%%%%%%%%%%%%%%%%%%%%%%%%%
%%%%%%%%%%%%%%%%%%%%%%%%%%%%%%%%%%%%%%%%%%%%%%%%%%%%%%%%%%%%%%%%%%%%%%%%%%%%%%%%%%
\subsection{Positivity and time monotonicity}\label{ptm}
%%%%%%%%%%%%%%%%%%%%%%%%%%%%%%%%%%%%%%%%%%%%%%%%%%%%%%%%%%%%%%%%%%%%%%%%%%%%%%%%%%
%%%%%%%%%%%%%%%%%%%%%%%%%%%%%%%%%%%%%%%%%%%%%%%%%%%%%%%%%%%%%%%%%%%%%%%%%%%%%%%%%%
A further property of $v$ is its time monotonicity which follows from the homogeneity of the operator $\Delta_\infty$ by a
result from B\'enilan \& Crandall \cite{BC}.
\begin{lem}\label{lem3.1}
  For $x \in \bar{\Omega}$, $s_1 \in \mathbb{R}$, $s_2 \in \R$ such that $s_1 \le s_2$, we have 
  $$v(s_1,x) \le v(s_2,x).$$
\end{lem} 
\proof Theorem~\ref{theo1.1} provides the well-posedness of \eqref{P} in $C_0 (\bar{\Omega})$ which is an ordered vector 
  space. As the comparison principle is valid for \eqref{P}-\eqref{PIC} by \cite[Theorem~2.3]{AS08} and the 
  infinity-Laplacian is homogeneous of degree 3, \cite[Theorem~2]{BC} implies 
  \begin{equation}\label{3.1.1}
    u(t+h,x) - u(t,x) \ge \left( \Big( \frac{t+h}{t} \Big)^{-1/2} -1 \right) u(t,x) \qquad\mbox{for } (t,x) \in (0,\infty)
    \times \bar{\Omega}, \, h>0. 
  \end{equation} 
  Hence, for any $(s,x) \in \mathbb{R} \times \bar{\Omega}$ and $h>0$, we obtain
  \begin{eqnarray*}
    v(s+h,x) - v(s,x) &=& e^{(s+h)/2} u(e^{s+h},x) - e^{s/2} u(e^s,x) \\
    &\ge& e^{(s+h)/2} \Big( \frac{e^{s+h}}{e^s} \Big)^{-1/2} u(e^s,x) - e^{s/2} u(e^s,x) = 0,
  \end{eqnarray*}
  which is the expected result.
\qed
The monotonicity of $v$ now enables us to prove that $v$ eventually becomes positive inside $\Omega$.
\begin{lem}\label{lem3.2}
  For any compact subset $K \subset \Omega$ there are $s_K >0$ and $\mu_K >0$ such that 
  \begin{equation}\label{3.2.1}
    v(s,x) \ge \mu_K >0 \qquad\mbox{in } [s_K, \infty) \times K.
  \end{equation}    
\end{lem}
\proof Three steps are needed to achieve the claimed result: we first prove that if $v(s, \cdot)$ is positive at one point
  of $\Omega$, then it becomes positive on a ``large'' ball centered around this point after a finite time. The second step
  is to prove that $v(s,\cdot)$ becomes eventually positive in $\Omega$ as $s \to \infty$, from which we deduce \eqref{3.2.1}
  in a third step.
 
  \medskip
  
  \textsl{Step 1:} Consider first $(t_0,x_0) \in (0,\infty) \times \Omega$ such that there are $\eps>0$ and $\delta>0$ with
  $B(x_0, \eps) \subset \Omega$ and
  \begin{equation}\label{3.2.3}
    u(t_0,x) \ge \delta >0 \qquad\mbox{for } x \in B(x_0,\eps).
  \end{equation} 
  Then, choosing $\alpha := \min \{ (4 \delta)^{1/3}, \eps^{2/3} \}$, $T:= (d(x_0,\partial \Omega)^6 / \alpha^9) -1 \ge 0$,
  and defining
  $$\mathcal{B}(t,x) := \frac{\alpha^3}{4}(t-t_0+1)^{-1/6} \left( 1- \alpha^{-2}|x-x_0|^{4/3} (t-t_0+1)^{-2/9} \right)_+^{3/2},
  \quad (t,x) \in [t_0,\infty) \times \mathbb{R}^N,$$
  we deduce from \cite[Proposition~1 and Corollary~1]{AJK09} that $\mathcal{B}$ is a viscosity solution to \eqref{P} in 
  $(t_0,t_0 +T) \times \Omega$. In addition, on the one hand, we have by \eqref{3.2.3}
  $$\mathcal{B}(t_0,x) \le \frac{\alpha^3}{4} \le \delta \le u(t_0,x) \qquad\mbox{for } x \in B(x_0,\eps)$$ 
  and
  $$\mathcal{B} (t_0,x) = 0 \le u(t_0,x) \qquad\mbox{for } x \in \bar{\Omega} \setminus B(x_0,\eps).$$
  On the other hand, we have $u(t,x) = \mathcal{B}(t,x) = 0$ for $(t,x) \in [t_0, t_0+T] \times \partial \Omega$ thanks to
  the choice of $T$, $\alpha$ and the properties of $\mathcal{B}$.
  The comparison principle \cite[Theorem~8.2]{CIL92} then implies $u \ge \mathcal{B}$ in $[t_0, t_0+T] \times \bar{\Omega}$.
  In particular, we have
  \begin{equation}\label{3.2.4}
    u(t_0+T,x) >0 \quad\mbox{for } x \in B(x_0, d(x_0,\partial \Omega)),
  \end{equation}
  where $T$ only depends on $\eps$ and $\delta$, but is independent of $x_0$ and $t_0$.
  
  \medskip
  
  \textsl{Step 2:} We next define the positivity set $\mathcal{P}(s)$ of $v(s, \cdot)$ for $s \ge 0$ by
  $$\mathcal{P}(s) := \{ x \in \Omega \, : \, v(s,x) >0 \}.$$
  Owing to the time monotonicity of $v$ (Lemma~\ref{lem3.1}), $(\mathcal{P}(s))_{s \ge 0}$ is a non-decreasing family of open
  subsets of $\Omega$ and
  $$\mathcal{P}_\infty := \bigcup_{s \ge 0} \mathcal{P}(s) \mbox{ is an open subset of } \Omega.$$
  Assume for contradiction that $\partial \mathcal{P}_\infty \cap \Omega \neq \emptyset$. Then there is $x_0 \in \partial 
  \mathcal{P}_\infty \cap \Omega$. Since $d(x_0, \partial \Omega) >0$ there is $y_0 \in \mathcal{P}_\infty$ such that $|y_0 - 
  x_0| \le d(x_0, \partial \Omega) /2 < d(y_0, \partial \Omega)$. Next, since $y_0 \in \mathcal{P}_\infty$, there is $s_0 >0$
  such that $v(s_0,y_0)>0$, that is $u(e^{s_0}, y_0) >0$. The previous step then guarantees the existence of $T \ge 0$, such
  that $u(e^{s_0}+T,x) >0$ for $x \in B(y_0, d(y_0, \partial \Omega))$. As $x_0 \in B(y_0, d(y_0, \partial \Omega))$, we 
  deduce from this that 
  $$v(\ln (e^{s_0} +T), x_0) = (e^{s_0}+T)^{1/2} u(e^{s_0}+T,x_0) >0,$$
  which contradicts the fact that $x_0 \in \partial \mathcal{P}_\infty$. Therefore, $\partial \mathcal{P}_\infty \cap \Omega
  = \emptyset$ and $\Omega$ is the union of the two disjoint open sets $\mathcal{P}_\infty$ and $\Omega \setminus \overline{
  \mathcal{P}_\infty}$. Since $\mathcal{P}_\infty \neq \emptyset$ by \eqref{1.1.1}, the connectedness of $\Omega$ implies
  \begin{equation}\label{3.2.5}
    \Omega = \mathcal{P}_\infty.
  \end{equation}  
  
  \medskip
  
  \textsl{Step 3:} Let $K$ be a compact subset of $\Omega$ and assume for contradiction that $K \not\subset \mathcal{P}(n)$
  for each $n \ge 1$. Then there is a sequence $(x_n)_{n \ge 1}$ in $K$ such that $v(n,x_n) = 0$ for $n \ge 1$ and we may
  assume without loss of generality that $x_n$ converges towards $x_\infty \in K$ as $n \to \infty$, thanks to the 
  compactness of $K$. Since $x_\infty \in \Omega$, it follows from \eqref{3.2.5} that there is $s_\infty>0$ such that 
  $v(s_\infty, x_\infty)>0$. Owing to the continuity of $v(s_\infty, \cdot)$ there are $\eps>0$ and $\delta>0$ such that
  $v(s_\infty,x) \ge \delta$ for $x \in B(x_\infty, \eps) \subset \Omega$. But then for $n$ large enough we have $n \ge
  s_\infty$ and $x_n \in B(x_\infty,\eps)$ and it follows from Lemma~\ref{lem3.1} and the previous bound that
  $$0 = v(n,x_n) \ge v(s_\infty,x_n) \ge \delta$$
  and a contradiction. Consequently, there is $n_K$ such that $K \subset \mathcal{P}(n_K)$ and 
  $$\mu_K := \min\limits_{x \in K} v(n_K,x) >0.$$ 
  Due to the time monotonicity of $v$, this implies \eqref{3.2.1}.   
\qed

%
%
%
%
%
%
%%%%%%%%%%%%%%%%%%%%%%%%%%%%%%%%%%%%%%%%%%%%%%%%%%%%%%%%%%%%%%%%%%%%%%%%%%%%%%%%%%
%%%%%%%%%%%%%%%%%%%%%%%%%%%%%%%%%%%%%%%%%%%%%%%%%%%%%%%%%%%%%%%%%%%%%%%%%%%%%%%%%%
\subsection{Convergence}\label{cv}
%%%%%%%%%%%%%%%%%%%%%%%%%%%%%%%%%%%%%%%%%%%%%%%%%%%%%%%%%%%%%%%%%%%%%%%%%%%%%%%%%%
%%%%%%%%%%%%%%%%%%%%%%%%%%%%%%%%%%%%%%%%%%%%%%%%%%%%%%%%%%%%%%%%%%%%%%%%%%%%%%%%%%
Having studied the positivity properties of $v$, we next turn to its behaviour near the boundary of $\Omega$ and first show the following lemma which is a modification of \cite[Lemma~10.1]{Li96}.

%%%%%%%%%%%%%%%%%%%%%%%%%%%%%%%%%%%%%%%%%%%%%%%%%%%%%%%%%%%%%%%%%%%%%%%%%%%%%%%%%%
\begin{lem}\label{led1} 
Consider $x_0\in\partial\Omega$, $\alpha\in (0,1/2)$, $\delta>0$, $B>0$, and define
$$
\psi_{\delta,B}(r) := \delta + B\ \left( r - \frac{r^2}{2} \right)\,, \quad r\in\R\,.
$$
Let $y_0\in \R^N$ be such that $|x_0-y_0|=R$ and $\Omega\cap B(y_0,R)=\emptyset$ (such a point $y_0$ exists according to the uniform exterior sphere condition \eqref{esc}). Introducing 
$$
U_{\alpha,x_0} := \left\{ x\in\Omega\ :\ R<|x-y_0|<R+\alpha \right\}
$$
and 
$$
w(s,x) := \psi_{\delta,B}\left( |x-y_0|-R \right)\,, \quad (s,x)\in [0,\infty)\times \overline{U_{\alpha,x_0}}\,,
$$
then $w$ is a supersolution to \eqref{v} in $(0,\infty)\times U_{\alpha,x_0}$ if $B\ge 2(1+\delta)$.
\end{lem}
%%%%%%%%%%%%%%%%%%%%%%%%%%%%%%%%%%%%%%%%%%%%%%%%%%%%%%%%%%%%%%%%%%%%%%%%%%%%%%%%%%

\proof To simplify notations, we set $\psi:=\psi_{\delta,B}$ and $U:=U_{\alpha,x_0}$. Since $\psi\in C^\infty(\R)$ and $y_0\not\in U$, the function $w$ is $C^\infty$-smooth in $(0,\infty)\times U$ and, if $(s,x)\in (0,\infty)\times U$, we have
\begin{equation}
\partial_s w(s,x) - \Delta_\infty w(s,x) - \frac{w(s,x)}{2} = - \left( {\psi'}^2\ \psi'' + \frac{\psi}{2} \right)\left( |x-y_0|-R \right).
\label{d3}
\end{equation}
Since $\alpha\in (0,1/2)$ and $B\ge 2$, we have for $r\in [0,\alpha]$
$$
- \left( {\psi'}^2\ \psi'' + \frac{\psi}{2} \right)(r) = B^3\ (1-r)^2 - \frac{B}{2}\ \left( r - \frac{r^2}{2} \right) - \frac{\delta}{2} \ge \frac{B^3}{8} - \frac{B}{4} - \frac{\delta}{2} \ge \frac{B-2\delta}{4}\,.
$$
Consequently, as $|x-y_0|-R\in [0,\alpha]$ for $(s,x)\in (0,\infty)\times U$, we deduce from \eqref{d3} and the above inequality that 
$$
\partial_s w(s,x) - \Delta_\infty w(s,x) - \frac{w(s,x)}{2}  \ge \frac{B-2\delta}{4} \ge 0\,,
$$ 
the last inequality following from the choice of $B$.
\qed

\medskip

As a consequence of Lemma~\ref{led1}, we have the following useful bound for $v$ on $\partial\Omega$.

%%%%%%%%%%%%%%%%%%%%%%%%%%%%%%%%%%%%%%%%%%%%%%%%%%%%%%%%%%%%%%%%%%%%%%%%%%%%%%%%%%
\begin{lem}\label{led2} 
Consider $\alpha\in (0,1/2)$ and define 
\begin{equation}
\label{d3b}
\omega(\alpha) := \sup{\{ v(0,x)\ :\ x\in\Omega \;\;\mbox{ and }\;\; d(x,\partial\Omega)<\alpha \}}\,.
\end{equation}
Then there is $\alpha_0\in (0,1/2)$ such that, for any $\alpha\in (0,\alpha_0)$ and $x_0\in\partial\Omega$, we have
\begin{equation}
\label{d4}
0 \le v(s,x) \le \omega(\alpha) + \frac{2C_1}{\alpha}\ |x-x_0|\,, \qquad (s,x)\in [0,\infty)\times (\bar{\Omega}\cap B(x_0,\alpha))\,,
\end{equation}
the constant $C_1$ being defined in \eqref{3.1a}.
\end{lem}
%%%%%%%%%%%%%%%%%%%%%%%%%%%%%%%%%%%%%%%%%%%%%%%%%%%%%%%%%%%%%%%%%%%%%%%%%%%%%%%%%%

\proof Consider $x_0\in\partial\Omega$ and let $y_0\in\R^N$ be such that $|x_0-y_0|=R$ and $\Omega\cap B(y_0,R)=\emptyset$, the existence of such a point $y_0$ being guaranteed by the uniform exterior sphere condition \eqref{esc}. With the notations of Lemma~\ref{led1}, we define
$$
w(s,x) := \psi_{\omega(\alpha),2C_1/\alpha}(|x-y_0|-R)\,, \quad (s,x)\in [0,\infty)\times \overline{U_{\alpha,x_0}}\,,
$$ 
the constant $C_1$ being defined in \eqref{3.1a} and observe that 
\begin{equation}\label{d4b}
B(x_0,\alpha)\cap \Omega \subset U_{\alpha,x_0} \subset \{ x\in\Omega\ :\ d(x,\partial\Omega)<\alpha \}\,.
\end{equation} 
On the one hand, it follows from \eqref{d3b} and \eqref{d4b} that
$$
w(0,x) \ge \omega(\alpha) \ge v(0,x)\,, \qquad x\in U_{\alpha,x_0}\,.
$$ 
On the other hand, if $(s,x)\in [0,\infty)\times\partial U_{\alpha,x_0}$, we have either $x\in\partial\Omega$ and $w(s,x)\ge 0=v(s,x)$ or $|x-y_0|=R+\alpha$ and
$$
w(s,x) = \psi_{\omega(\alpha),2C_1/\alpha}(\alpha) \ge \frac{2C_1}{\alpha}\ \left( \alpha - \frac{\alpha^2}{2} \right) \ge C_1\ge v(s,x)
$$
by \eqref{3.1a}. Furthermore, since $v(0,x)=0$ on $\partial\Omega$, $\omega(\alpha)$ converges to $0$ as $\alpha \searrow 0$ and there is thus $\alpha_0\in (0,1/2)$ such that $2C_1/\alpha\ge 2(1+\omega(\alpha))$ for $\alpha\in (0,\alpha_0)$. This condition implies that $w$ is a supersolution to \eqref{v} in $(0,\infty)\times U_{\alpha,x_0}$ by Lemma~\ref{led1}. According to the above analysis, we are in a position to apply the comparison principle \cite[Theorem~8.2]{CIL92} to conclude that 
$$
v(s,x)\le w(s,x)\,, \qquad (s,x)\in [0,\infty)\times \overline{U_{\alpha,x_0}}\,.
$$
In particular, if $(s,x)\in [0,\infty)\times (\bar{\Omega}\cap B(x_0,\alpha))$, the above inequality, \eqref{d4b}, and the properties of $y_0$ entail that
\begin{eqnarray*}
v(s,x) & \le & \omega(\alpha) + \frac{2C_1}{\alpha}\ (|x-y_0|-R) \\
& \le & \omega(\alpha) + \frac{2C_1}{\alpha}\ (|x-x_0|+|x_0-y_0|-R) \\
& \le & \omega(\alpha) + \frac{2C_1}{\alpha}\ |x-x_0|\,,
\end{eqnarray*}
whence \eqref{d4}. 
\qed

\medskip

\proofc{of Theorem~\ref{theo1.2}}
For $\eps\in (0,1)$, we define
$$
V_\eps(s,x) := v\left( \frac{s}{\eps} , x \right)\,, \qquad  (s,x)\in [0,\infty)\times\bar{\Omega}\,,
$$
and the half-relaxed limits
$$
V_*(x) := \liminf_{(\sigma,y,\eps)\to (s,x,0)} V_\eps(\sigma,y)\,, \qquad V^*(x) := \limsup_{(\sigma,y,\eps)\to (s,x,0)} V_\eps(\sigma,y)
$$
for  $(s,x)\in (0,\infty)\times\bar{\Omega}$. These functions are well-defined by \eqref{3.1a}, indeed do not depend on $s>0$, and the stability result for (discontinuous) viscosity solutions ensures that 
\begin{eqnarray}
& & V_* \;\;\mbox{ is a supersolution to }\;\; -\Delta_\infty z - \frac{z}{2} =0 \;\;\mbox{ in }\;\; \Omega\,, \label{d7} \\
& & V^* \;\;\mbox{ is a subsolution to }\;\; -\Delta_\infty z - \frac{z}{2} =0 \;\;\mbox{ in }\;\; \Omega\,. \label{d8}
\end{eqnarray} 
In addition, it follows from \eqref{3.1a} and \eqref{d4} that
\begin{equation}
0 \le V_*(x) \le V^*(x) \le C_1\,, \quad x\in\bar{\Omega}\,, \label{d5}
\end{equation}
and, for all $(x_0,\alpha) \in \partial \Omega \times (0,\alpha_0)$,
\begin{equation}
0 \le V_*(x) \le V^*(x) \le \omega(\alpha) + \frac{2C_1}{\alpha}\ |x-x_0|\,, \quad x\in \bar{\Omega}\cap B(x_0,\alpha)\,. \label{d6}
\end{equation}
In particular, \eqref{d6} guarantees that $0\le V_*(x_0)\le V^*(x_0) \le \omega(\alpha)$ for all $x_0\in\partial\Omega$ and $\alpha\in (0,\alpha_0)$. Since $\omega(\alpha)\to 0$ as $\alpha \searrow 0$, we end up with
\begin{equation}\label{d6b}
V_*(x)=V^*(x) = 0\,, \qquad x\in\partial\Omega\,.
\end{equation}
We finally infer from Lemma~\ref{lem3.2} that
\begin{equation}\label{d9}
  V_* (x) >0 \qquad\mbox{for } x \in \Omega.
\end{equation}
We are then in the position to apply Lemma~\ref{lem2.1} to obtain that $V^* \le V_*$. Recalling \eqref{d7}, 
\eqref{d8}, \eqref{d5}, and \eqref{d6b} we conclude that $V_* = V^* \in C_0 (\bar{\Omega})$ is a viscosity solution to 
$-\Delta_\infty z - z/2 = 0$ in $\Omega$. We have thus proved that $f_\infty := V^*$ is a positive viscosity solution to
\eqref{1.2.2} and it is the only one by Corollary~\ref{cor2.2}. In addition, it follows from the identity $V^* = V_* = 
f_\infty$ and \cite[Lemme~4.1]{Bl94} (see also \cite[Lemma~5.1.9]{BdCD97}) that
$$\lim\limits_{\eps \searrow 0} \| V_\eps (2) - f_\infty \|_{L^\infty (\Omega)} = 0.$$  
In other words,
\begin{equation}\label{d100}
  \lim\limits_{s \to \infty} \| v(s) - f_\infty\|_{L^\infty (\Omega)} = 0,
\end{equation}  
which is equivalent to \eqref{1.2.1} by \eqref{3.1}.
\qed
\proofc{of Corollary~\ref{cor1.3}} The claim now follows from Theorem~\ref{theo1.2} and Lemma~\ref{lem3.1}. Indeed, we have
$v(s,\cdot) \le v(\sigma, \cdot)$ for $-\infty < s \le \sigma < \infty$. Letting $\sigma \to \infty$ and using \eqref{d100}
lead us to $v(s, \cdot) \le f_\infty$ for any $s \in \mathbb{R}$, which is nothing but \eqref{1.3.1} once written in terms
of $u$.  
\qed
%
%
%
%
%
%
%%%%%%%%%%%%%%%%%%%%%%%%%%%%%%%%%%%%%%%%%%%%%%%%%%%%%%%%%%%%%%%%%%%%%%%%%%%%%%%%%%
%%%%%%%%%%%%%%%%%%%%%%%%%%%%%%%%%%%%%%%%%%%%%%%%%%%%%%%%%%%%%%%%%%%%%%%%%%%%%%%%%%
\mysection{Additional positivity properties}\label{app}
%%%%%%%%%%%%%%%%%%%%%%%%%%%%%%%%%%%%%%%%%%%%%%%%%%%%%%%%%%%%%%%%%%%%%%%%%%%%%%%%%%
%%%%%%%%%%%%%%%%%%%%%%%%%%%%%%%%%%%%%%%%%%%%%%%%%%%%%%%%%%%%%%%%%%%%%%%%%%%%%%%%%%
First we state an extension of Lemma~\ref{lem3.2} which shows that $u$ is indeed positive in $\Omega$ after a finite time
provided that $\Omega$ additionally satisfies a uniform interior sphere condition:
\begin{equation}\label{4.1}
  \begin{split} & \mbox{There is } R_0>0 \mbox{ such that for any } x_0 \in \partial \Omega \mbox{ there is } y_0 \in \Omega  
  \\
  & \mbox{such that } |y_0-x_0| = R_0 \mbox{ and } B(y_0,R_0) \subset \Omega.
  \end{split}
\end{equation}
\begin{lem}\label{lem3.2a}
  Assume that $\Omega \subset \mathbb{R}^N$ is a bounded domain satisfying \eqref{esc} and \eqref{4.1} and that $u_0$ 
  fulfills \eqref{ic}. If $u$ denotes the viscosity solution to \eqref{P}-\eqref{PIC}, 
  then there is $t_1 \in (0,\infty)$ such that 
  \begin{equation}\label{3.2.2}
    u(t,x) >0 \qquad\mbox{in } [t_1, \infty) \times \Omega.
  \end{equation}  
\end{lem}
\proof Let $v$ be defined by \eqref{3.1} and set 
  $$K := \Big\{ x \in \Omega \, : \, d(x, \partial \Omega) \ge \frac{R_0}{2} \Big\} \quad\mbox{ and }\quad 
  M:= \{ x \in \Omega \, : \, d(x, \partial \Omega) = R_0\}.$$ 
  Since $K$ is a compact subset of $\Omega$, we have
  \begin{equation}\label{3.2.1a}
    v(s,x) \ge \mu_K >0 \qquad\mbox{in } [s_K, \infty) \times K
  \end{equation}  
  for some $s_K >0$ and $\mu_K>0$ by Lemma~\ref{lem3.2}. Thus, setting $t_0 := e^{s_K}$, 
  $\eps := R_0 / 2$ and $\delta := t_0^{-1/2} \mu_K$, \eqref{3.2.3} is valid for any $x_0 \in M$. Then the first step
  of the proof of Lemma~\ref{lem3.2} implies the existence of $T>0$ which is independent of $x_0 \in M$ such that 
  \eqref{3.2.4} is fulfilled for any $x_0 \in M$. Thus, we conclude that
  $$v(s_0,x) >0 \quad\mbox{for } x \in \tilde{M} := \bigcup_{x_0 \in M} B(x_0,R_0),$$       
  where $s_0 := \ln (t_0 +T) > s_K$. As \eqref{4.1} implies $\tilde{M} \cup K = \Omega$ (see e.g. \cite[Section~14.6]{GT01}),
  we deduce from Lemma~3.1 and \eqref{3.2.1a} that
  $$v(s,x) >0 \qquad\mbox{in } [s_0, \infty) \times \Omega.$$
  By \eqref{3.1}, this shows \eqref{3.2.2} with $t_1 := e^{s_0}$. 
\qed
Having shown that $u$ is positive in $\Omega$ after a finite or infinite time, we next show that the expansion of the
positivity set of $u(t, \cdot)$ may take some time to be initiated. 
\begin{prop}\label{pre1}
Consider $u_0\in C_0(\bar{\Omega})$ and define its positivity set $\mathcal{P}_0$ by 
$$
\mathcal{P}_0 := \{ x \in \Omega\ : \ u_0(x)>0 \}\,.
$$
If $x_0\in \Omega \cap \partial \mathcal{P}_0$ is such that 
\begin{equation}\label{e1}
u_0(x) \le a\ |x-x_0|^2\,, \qquad x\in B(x_0,\delta)\subset \Omega\,,
\end{equation}
for some $\delta>0$ and $a>0$, then there is $\tau(x_0)>0$ such that $u(t,x_0)=0$ for $t\in [0,\tau(x_0))$.
\end{prop}
In other words, the so-called waiting time 
$$\tau_w (x_0) := \inf \{ t>0 \, : \, u(t,x_0) >0 \}$$
of $u$ at $x_0 \in \Omega$ is positive if $u_0$ satisfies \eqref{e1}. In addition, it is finite by Lemma~\ref{lem3.2}.
This waiting time phenomenon is typical for degenerate parabolic equations, see \cite{DPGG03, Va07} and the references
therein.

The proof of Proposition~\ref{pre1} relies on the construction of supersolutions as in \cite[Theorem~8.2]{Kn77} which we describe now.

\begin{lem}\label{lee2}
Consider $x_0\in\Omega$ and $T>0$ and define
$$
S_T(t,x) := \frac{|x-x_0|^2}{4 (T-t)^{1/2}}\,, \qquad (t,x)\in [0,T)\times\bar{\Omega}\,.
$$
Then $S_T$ is a supersolution to \eqref{P} in $(0,T)\times\Omega$.
\end{lem}

\proof We first note that $S_T\in C^2([0,T)\times\bar{\Omega})$. For $(t,x)\in (0,T)\times\Omega$, we compute
\begin{eqnarray*}
\partial_t S(t,x) - \Delta_\infty S(t,x) & = & \frac{|x-x_0|^2}{8 (T-t)^{3/2}} - \frac{\langle x-x_0 , x-x_0 \rangle}{8 (T-t)^{3/2}} = 0 
\end{eqnarray*}
and readily obtain the expected result. \qed

\medskip

\proofc{of Proposition~\ref{pre1}} Define
$$
T := \min{\left\{ \frac{1}{16a^2} , \frac{\delta^4}{16 C_1^2} \right\}}\,.
$$
According to Lemma~\ref{lee2}, the function $S_T$ is a supersolution to \eqref{P} in $(0,T)\times B(x_0,\delta)$. In addition, the choice of $T$ and \eqref{e1} guarantee that
$$
S_T(0,x) = \frac{|x-x_0|^2}{4 T^{1/2}} \ge a\ |x-x_0|^2 \ge u_0(x)\,, \quad x\in B(x_0,\delta)\,,
$$
while we infer from the choice of $T$ and \eqref{1.1.1} that, for $(t,x)\in (0,T)\times \partial B(x_0,\delta)$
$$
S_T(t,x) = \frac{\delta^2}{4 (T-t)^{1/2}} \ge \frac{\delta^2}{4 T^{1/2}} \ge C_1 \ge u(t,x) \,.
$$
The comparison principle \cite[Theorem~8.2]{CIL92} then entails that $S_T(t,x)\ge u(t,x)$ for $(t,x)\in [0,T)\times B(x_0,\delta)$. In particular, $0\le u(t,x_0)\le S_T(t,x_0)=0$ for $t\in [0,T)$, and the proof of Proposition~\ref{pre1} is complete. \qed

%%%%%%%%%%%%%%%%%%%%%%%%%%%%%%%%%%%%%%%%%%%%%%%%%%%%%%%%%%%%%%%%%%%%%%%%%%%%%%%%%%
%%%%%%%%%%%%%%%%%%%%%%%%%%%%%%%%%%%%%%%%%%%%%%%%%%%%%%%%%%%%%%%%%%%%%%%%%%%%%%%%%%
\bigskip

\textbf{Acknowledgements}

This work was done while C.~Stinner held a one month invited position at the Institut de Math\'{e}-matiques de
Toulouse, Universit\'{e} Paul Sabatier - Toulouse III. He would like to express his gratitude for the invitation, support, and 
hospitality.  
 
%%%%%%%%%%%%%%%%%%%%%%%%%%%%%%%%%%%%%%%%%%%%%%%%%%%%%%%%%%%%%%%%%%%%%%%%%%%%%%%%%%
%%%%%%%%%%%%%%%%%%%%%%%%%%%%%%%%%%%%%%%%%%%%%%%%%%%%%%%%%%%%%%%%%%%%%%%%%%%%%%%%%%
%
%
%
%%%%%%%%%%%%%%%%%%%%%%%%%%%%%%%%%%%%%%%%%%%%%%%%%%%%%%%%%%%%%%%%%%%%%%%%%%%%%%%%%%
%%%%%%%%%%%%%%%%%%%%%%%%%%%%%%%%%%%%%%%%%%%%%%%%%%%%%%%%%%%%%%%%%%%%%%%%%%%%%%%%%%

%
%
%
%
\end{document}